\documentclass[11pt]{article}
\usepackage{amsmath,latexsym, color, amsbsy,amssymb}
\usepackage{graphicx}
\usepackage{latexsym}
\def\forall{\hbox{for all}~}

\def\ve{\varepsilon}
\def\n{\noindent}
\def\ds{\displaystyle}

\def\R{\mathbb{R}}

\def\v{\vskip 1em}
\def\O{{\cal O}}
\def\C{{\cal C}}

\def\begi{\begin{itemize}}
	\def\endi{\end{itemize}}
\def\bega{\begin{array}}
	\def\enda{\end{array}}

\def\Hat{\widehat}

\def\bel{\begin{equation}\label}
	\def\eeq{\end{equation}}
\def\sqr#1#2{\vbox{\hrule height .#2pt
		\hbox{\vrule width .#2pt height #1pt \kern #1pt
			\vrule width .#2pt}\hrule height .#2pt }}
\def\square{\sqr74}
\def\endproof{\hphantom{MM}\hfill\llap{$\square$}\goodbreak}

\newtheorem{theorem}{Theorem}[section]
\newtheorem{corollary}[theorem]{Corollary}
\newtheorem{definition}[theorem]{Definition}
\newtheorem{remark}[theorem]{Remark}

\newtheorem{lemma}[theorem]{Lemma}

\begin{document}
	\title{\bf   Two-Point Boundary Value Problems for Quasi-Monotone Dynamical Systems}\v

	\author{\it Lorena Bociu, Khai T. Nguyen, and Madhumita Roy\\
		\\
		{\small Department of Mathematics, North Carolina State University}\\
		{\small e-mails: ~lvbociu@ncsu.edu, ~khai@math.ncsu.edu,  ~mroy5@ncsu.edu }
	}

	\maketitle
	\begin{abstract} This paper studies the existence of minimal solutions to two-point boundary value problems for quasi-monotone dynamical systems. Specifically, the pointwise infimum of all supersolutions is shown to coincide with the minimal solution. This result is then applied to establish a non-uniqueness result for strong stable solutions to a class of mean field games with a continuum of players.
		\quad\\
		\quad\\
		{\footnotesize
			\medskip
			
			\n {\bf AMS Mathematics Subject Classification.} 34C12, 34D20, 49N80, 	49L12, 91A16
		}
	\end{abstract}

	\section{Introduction}\label{sec:1}
	\setcounter{equation}{0}
	We consider a two-point boundary value problem for the system of ordinary differential equations (ODEs)
	\bel{1}\left\{\bega{rl}\dot x &= f(t,x,y),\\[2mm]
	\dot y&= g(t,x,y),\enda\right.\qquad a.e.~t\in [0,T]\,,\eeq
	with boundary conditions given by
	\bel{2}x(0)~=~\bar x\,,\qquad\qquad y(T)~=~\bar y\,,\eeq
	where $x: [0, T] \to \R^{N_1}$,  $x(t) =(x_1(t),\ldots,x_{N_1}(t))$ and $y: [0, T] \to \R^{N_2}$,  $y(t)=(y_1(t),\ldots,y_{N_2}(t))$ are absolutely continuous functions in $[0,T]$ for \({N_1, N_2 \in \mathbb{N}}\). 
	Here,  we shall assume the following standard hypothesis: 
	\begin{itemize}
		\item [{\bf (A).}]  The functions $f=(f_1,\cdots,f_{N_1}):\R\times \R^{N_1}\times\R^{N_2}\to \R^{N_1}$ and $g=(g_1,\cdots, g_{N_2}):\R\times\R^{N_1}\times\R^{N_2}\to\R^{N_2}$ are  continuous w.r.t $x$ and $y$, measurable w.r.t $t$, and locally bounded, i.e., for all $M>0$, it holds
		\bel{b}
		\sup_{t\in [0,T], |x|,|y|\leq M} \max\{|f(t,x,y)|, |g(t,x,y)|\}~<~+\infty.
		\eeq
	\end{itemize}
	
	The existence of Carath\'eodory solutions theory for two-point boundary value problems has been extensively studied by various authors across diverse cases arising from different applications.
	\begin{definition}[Carath\'eodory Solution]
		The absolutely continuous map $t\mapsto (x(t), y(t))$
		is a Carath\'eodory solution to (\ref{1})-(\ref{2}) if  for all $t\in [0,T]$, the following identities hold:
		\[
		x(t)~=~\bar{x}+\int_{0}^{t}f(s,x(s),y(s))ds,\qquad y(t)~=~\bar{y} - \int_{t}^{T}g(s,x(s),y(s)) ds.
		\] 
	\end{definition}
	In the classical work of Bucy \cite{RS}, the focus is on linear Hamiltonian systems subject to two-point boundary conditions. Bucy proves the uniqueness of solution under certain positive semidefinite assumption.  Kelevedjiev~\cite{K} proves existence results for nonlinear two-point boundary value problems of the form \(\ddot{x} = f(t,x,\dot{x}),\) by combining a priori bounds on both \(x\) and \(\dot{x}\) with topological transversality arguments. Unlike methods relying on global Lipschitz conditions, this approach works under weak regularity and nonstandard growth. Fonda, Sfecci, and Toader \cite{A} extend Picard’s 1893 upper–lower solution theory to planar nonlinear systems with boundary data given by two lines in the plane rather than fixed points, accommodating second-order equations. Using degree theory, they prove the existence of a solution between ordered upper and lower bounds, with techniques particularly suited to problems exploiting monotonicity. Later, Fonda and Ortega \cite{A1} treated two-point boundary value problems using variational and symplectic techniques to establish the geometric multiplicity of solutions for Hamiltonian systems on a cylinder. Their framework shows that, for Hamiltonian problems with geometric structure, one can exploit critical-point theory or Poincaré-map arguments to produce multiple geometrically distinct solutions. %

	From a different aspect, the system of ODEs (\ref{1})  represents  the forward-backward system of characteristics for nonconservative systems of transport of the form
	\bel{Tr}
	\partial_tu  + f(t,x,u)\cdot\nabla u~ =~g(t,x,u),
	\eeq
	which are important in the study of mean field games (MFG). These are models for large populations of interacting rational agents, which strategize in order to optimize an outcome, based on the collective behavior of the remaining population, while subject to environmental influences.
	More precisely, (\ref{Tr}) is the master equation of mean field games with a (finite) discrete state space with $N_1=N_2$.  In order to draw the analogy to MFG PDE systems and the master equation in a continuum state space  it is convenient to represent the characteristics as the system of ODEs (\ref{1}) \cite{C}. It is known that if $f$ and $g$ are smooth and $T$ is sufficiently small then (\ref{1})-(\ref{2}) has a unique Carath\'eodory solution \cite{LS}. However, the argument fails for arbitrarily long time intervals, in view of the coupling between $x$ and the terminal condition on $y$. In the case ${N_1}={N_2}$, if the maps   $(g,f):\R^{N_1}\times\R^{N_1} \to \R^{N_1}\times\R^{N_1}$ are smooth and monotone, with at least one component strictly monotone, then   the system (\ref{1})-(\ref{2}) again admits a unique solution on any time interval, provided the terminal condition $\bar{y}=\psi(\bar{x})$  is given by a smooth and monotone function $\psi$  \cite{BL}.  This setting is directly analogous to the monotonicity condition introduced by Lasry and Lions for mean field game systems with a continuum state space in \cite{BC, L}.

	In this paper, inspired by a study on the existence of multiple strong solutions to a class of mean field games with a continuum of players in \cite{BN}, where each player's state evolves according to a controlled ODE, we   introduce the concept of a minimal solution to the two-point boundary value problem (\ref{1})-(\ref{2}). On the spaces $\R^{N_1}$ and $\R^{N_2}$, we consider the partial orderings given by
	$$\bega{rl}x\preceq\tilde x\qquad &\hbox{if and only if}\qquad x_i\leq\tilde x_i
	\quad\forall i=1,\ldots,N_1,\\[1mm]
	y\preceq\tilde y\qquad &\hbox{if and only if}\qquad y_j\leq\tilde y_j
	\quad\forall j=1,\ldots,N_2.\enda$$

	\begin{definition}[Minimal Solution]
		A Carathéodory solution \( (x(t), y(t)) \) of \eqref{1}–\eqref{2} is said to be a \textit{minimal solution} if, for any other Carathéodory solution \( (x_0(t), y_0(t)) \) of \eqref{1}–\eqref{2}, the following holds for all \( t \in [0, T] \):
		\[
		x(t) \preceq x_0(t) \quad \text{and} \quad y(t) \preceq y_0(t).
		\]
	\end{definition}
	Throughout the paper, we work under the following assumptions on the monotonicity  for the nonlinearities $f$ and $g$: 
	\begin{itemize}
		\item [{\bf (M).}] The system of ODEs (\ref{1}) satisfies  the {\it quasi-monotonicity property}:
		\begin{itemize}
			\item [{\bf (M$_f$).}]  For every $i = \overline{1, N_1}$, $f_i$ is non-decreasing  w.r.t $x_k$ and $y_j$,  for $j =\overline{1, N_2}$ and $k\neq i$.
			\item [{\bf (M$_g$).}] For every $j =\overline{1, N_2} $, $g_j$  is non-increasing w.r.t $x_i$ and $y_k$,  for $i = \overline{1, N_1}$ and $k\neq j$.
		\end{itemize}
	\end{itemize}
	To establish a well-posedness result for the minimal solution, we use the classical concept of a supersolution, introduced by Picard in 1983 for boundary value problems and then extended to various ODEs (see, e.g.,   \cite{B, A, N}, and PDEs (see, e.g., \cite{B, Evans, LS}).
	
	\begin{definition}[Supersolution]\label{dfnss}
		The absolutely continuous map $t\mapsto (x(t), y(t))$
		is a supersolution to (\ref{1})-(\ref{2}) 
		if the following conditions hold:
		\begi
		\item[(i)] For a.e.~$t\in [0,T]$ the following inequalities are satisfied:
		\bel{supers}\left\{\bega{rl}\dot x &\succeq~ f(t,x,y),\\[2mm]
		\dot y&\preceq ~g(t,x,y).\enda\right.\eeq
		\item[(ii)] The boundary data satisfy $\bar x\preceq x(0)$ and $ \bar y
		\preceq y(T)$.
		\endi
	\end{definition}

	{\bf Summary of Main Results.} 
	In Section 2, using the classical  supersolution approach, we prove that the two-point boundary value problem (\ref{1})-(\ref{2}) satisfying the quasi-monotonicity property {\bf (M)} admits a minimal solution, provided that the set of supersolutions is nonempty and uniformly bounded from below. Then we construct an example demonstrating that the boundedness assumption on the set of supersolutions is essential. This naturally leads to the question of identifying structural conditions on the underlying dynamics that guarantee such uniform boundedness of the supersolution set. Building on these observations, we  establish a general sufficient condition for the existence of a minimal solution to (\ref{1})-(\ref{2}). We remark that the minimal solution is  unique, but  multiple Carath\'eodory solutions  may exist under the same data.

	In Section 3, we consider a class of deterministic mean field games involving a continuum of players, where the state of each player evolves according to a simple linear controlled ordinary differential equation. Applying our results from Section 2 to the two-point boundary value problem derived from the Pontryagin Maximum Principle, we demonstrate that the corresponding mean field game admits at least two distinct solutions: one that is stable in the classical sense, and another one that is a minimal solution and asymptotically stable. This multiplicity of stable equilibria highlights the sensitivity of the game dynamics to initial configurations and control costs, and illustrates how nonuniqueness can persist even in structurally simple mean field models.

	\section{Existence  of the minimal solutions}
	In this section we  establish existence of minimal Carath\'eodory solutions for the two-point boundary value problem (\ref{1})-(\ref{2}) satisfying the  quasi-monotonicity property {\bf (M)} and {\bf (A)}. Towards this goal, we first prove the following lemma: 
	\begin{lemma}
		Assume that the system of ODEs (\ref{1})-(\ref{2}) satisfies {\bf (A)}. If  $(x^\flat, y^\flat)$ and $(x^\sharp, y^\sharp)$ are two supersolutions then the pointwise minimum $(x,y)$ defined as
		\bel{min}\bega{rl}
		x_i(t)&\doteq~\min~\bigl\{x_i^\flat (t),\, x_i^\sharp(t)\bigr\},
		\qquad i=1,\ldots,N_1,\\[4mm]
		y_j(t)&\doteq~\min~\bigl\{y_j^\flat (t),\, y_j^\sharp(t)\bigr\},
		\qquad j=1,\ldots,N_2,\enda\eeq
		is also a supersolution.
	\end{lemma}
	{\bf Proof.} By definition (\ref{min}), both $x=(x_1,x_2,...,x_{N_1})$ and $y=(y_1,y_2,...,y_{N_2})$ are absolutely continuous and satisfy the boundary conditions requirements from  Definition \ref{dfnss}, part (ii).  Hence, we only need to check  that they satisfy the inequalities in Definition \ref{dfnss}, part (i). Fix any $i\in\{1,...,{N_1}\}$ and let $\bar{t}\in [0,T]$ be a differentiable point of $x_i,x^{\flat}_i$ and $x_i^{\sharp}$. Without loss of generality, assuming $x_i(\bar{t})=x_i^{\flat}(\bar{t})$, it follows that 
	\[
	\dot{x}_i(\bar{t})~=~\dot{x}^{\flat}_i(\bar{t})~\geq~f_i(x^{\flat}(\bar{t}),y^{\flat}(\bar{t}))\,.
	\]
	Given that  $x_k(\bar{t})~\leq~x^{\flat}_k(\bar{t})$ for all $k\in\{1,2,...,N_1\}\backslash i$ and $y_j(\bar{t})\leq y^{\flat}_j(\bar{t})$ for all $j\in\{1,2,...,N_2\}$, using the quasi-monotonicity properties of $f$ in {\bf (M$_f$)}, we derive 
	\[
	\dot{x}_i(\bar{t})~\geq~f_i(\bar{t},x^{\flat}(\bar{t}),y^{\flat}(\bar{t}))~\geq~f_i(\bar{t},x(\bar{t}),y(\bar{t})).
	\]
	Since this inequality holds for all $i\in \{1,...,N_1\}$ and for almost every $\bar{t}\in [0,T]$, it holds that
	\[
	\dot{x}(t)~\succeq~f(t,x(t),y(t)),\qquad\mathrm{for}\ a.e.\ t\in [0,T]\,.
	\]
	Similarly, by using the quasi-monotonicity of $g$ in {\bf (M$_g$)}, one can show that 
	\[
	\dot{y}(t)~\preceq~g(t,x(t),y(t)),\qquad\mathrm{for}\ a.e.\ t\in [0,T],
	\]
	which completes the  proof of the lemma. 
	\endproof
	\medskip

	Our main result is stated below. 
	\begin{theorem}\label{T1}
		Consider the boundary value problem (\ref{1})-(\ref{2}) under the  quasi-monotonicity property {\bf (M)} and {\bf (A)}. Assume that the family of supersolutions is nonempty and  uniformly bounded below. Then the pointwise infimum of all supersolutions is a minimal  solution to (\ref{1})-(\ref{2}).
	\end{theorem}
	{\bf Proof.} Assume that  $(x^*(t),y^*(t))$ is the pointwise infimum of all supersolutions of the two-point boundary value problem(\ref{1})-(\ref{2}).

	In the first two steps of the proof, we construct a decreasing sequence of supersolutions $\big(x^n(t),y^n(t)\big)_{n\geq 1}$ of (\ref{1})-(\ref{2}) such that $\forall i=1,\ldots, N_1$ and $j=1,\ldots, N_2,
	$
	\bel{inf-1}
	x_i^{*}(t)~=~\lim_{n\to\infty} x_i^{n}(t),\qquad\qquad y_j^{*}(t)~=~\lim_{n\to\infty} y_j^{n}(t),\qquad\forall t\in [0,T]\,.
	\eeq
	{\bf 1.} Fix a time $\tau\in [0,T]$ and  let $(x^n,y^n)$ be a sequence of supersolutions of (\ref{1})-(\ref{2}) such that $\forall i=1,\ldots, N_1$ and $j=1,\ldots, N_2,
	$
	\[
	\lim_{n\to\infty}x_i^n(\tau)~=~x_i^*(\tau),\qquad \lim_{n\to\infty}y_j^n(\tau)~=~y_j^*(\tau).
	\]
	From Lemma 1, by replacing each couple functions $\big(x^n,y^n\big)$ with $\big(\tilde{x}^n,\tilde{y}^n\big)$ which is defined $\forall i=1,\ldots, N_1$ and $j=1,\ldots, N_2,
	$ by
	\[
	\tilde{x}^n_i(t)~=~\min_{k\in\{1,2,...,n\}}x^k_i(t),\qquad \tilde{y}^n_j(t)~=~\min_{k\in\{1,2,...,n\}}y^k_j(t),\quad\forall t\in [0,T],\]
	we can also assume that both $(x^n)_{n \geq 1}$ and $(y^n)_{n\geq 1}$ are  decreasing sequences. In this case,  since family of supersolutions of (\ref{1})-(\ref{2}) is uniformly bounded below, there is a constant $M>0$ such that
	\bel{b-xn}
	\sup_{t\in [0,T]} \max\{|x^n(t)|, \big|y^n(t)\big|\}~\leq~M.
	\eeq
	Using Assumption (A), we set 
	$$\ds M_1~\doteq~\sup_{t\in [0,T],|x|,|y|\leq M} \max\{|f(t,x,y)|,|g(t,x,y)|\}$$ 
	so that
	\[
	\dot{x}_i^n(t)~\geq~-M_1\qquad\mathrm{and}\qquad\dot{y}^n_j(t)~\leq~M_1\qquad a.e.\ t\in [0,T],
	\]
	In particular, this yields  
	\bel{One-side}
	\begin{cases}
		x^n_i(t)-x^n_i(\tau)~\geq~-M_1(t-\tau),\\[2mm]
		y^n_j(t)-y^n_j(\tau)~\leq~M_1(t-\tau),
	\end{cases}
	\qquad 0\leq \tau<t\leq T.
	\eeq
	Thus, for all $i\in \{1,2,...,N_1\}$ and $j\in\{1,2,...,N_2\}$, it holds 
	\bel{onesideL-x*}
	\begin{split}
		x_i^*(\tau)-x_i^*(t)~\geq~\liminf_{n\to\infty} \big(x^n_i(\tau)-x^n_i(t)\big)~\geq~M_1\cdot (t-\tau),\qquad t\in [0,\tau],
	\end{split}
	\eeq
	and
	\bel{onesideL-y*}
	y_j^*(t)-y_j^*(\tau)~\leq~\limsup_{n\to\infty} \big(y^n_j(t)-y^n_j(\tau)\big)~\leq~M_1\cdot (t-\tau),\qquad t\in [\tau,T].
	\eeq
	
	Since the above is true for every $\tau\in [0,T]$,  both $t\mapsto M_1t+x_i^*(t)$ and $t\mapsto M_1t-y_j^*(t)$ are increasing in the interval $[0,T]$. As a consequence, the set $\mathcal{D}$ of discontinuities of $x^*$ and $y^*$ in $[0,T]$ is at most countable. 
	\medskip

	{\bf  2.} From  Step 1, by a standard argument, we can construct  a decreasing sequence of supersolutions $\big(x^n,y^n\big)_{n\geq 1}$ which converges to $(x^*,y^*)$ at every rational time and $T$, and at every time  where  $x^*,y^*$ are discontinuous.  Indeed, a decreasing sequence on the countable set $\mathbb{Q}$ could be constructed in the following way: choose $t_1 \in \mathbb{Q}$ and a decreasing sequence $\{x_i^m(t_1)\}$ for some fixed $i$. Note that for another $t_2 \in \mathbb{Q}$ the sequence $\{x_i^m(t_2)\}$ might not be decreasing. Therefore next we construct a subsequence $\{x_i^{m_k}\} \subset \{x_i^m\}$ which will be decreasing on $t_2$ and so on.

	From the definition of $x^*,y^*$,  for all $i\in\{1,2,...,N_1\}$ and $j\in\{1,2,...,N_2\}$, one has
	\[
	x_i^*(\tau)~\leq~\lim_{n\to\infty}x^n_i(\tau)\qquad\mathrm{and}\qquad y_j^*(\tau)~\leq~\lim_{n\to\infty}y^n_j(\tau)\qquad\forall \tau\in [0,T].
	\]
	Now let $\tau\in (0,T)\backslash\mathcal{D}$.  Assume that  $x^*_{i}(\tau)<\ds\lim_{n\to\infty}x^n_i(\tau)-\ve$ for some $\ve>0$, $i\in\{1,2,...,N_1\}$. Since $x^*$ is continuous at $\tau$, we can choose  a rational time $\tau<t<T$ such that 
	\[
	x_i^*(t)~<~x_{i}^*(\tau)+{\ve\over 3}\qquad\mathrm{and}\qquad M_1(t-\tau)~<~{\ve\over 3}.
	\]
	Recalling (\ref{One-side}), we derive
	\[
	\begin{split}
		x^*_{i}(\tau)&~\leq~\lim_{n\to\infty}x_i^n(\tau)-\ve~\leq~\lim_{n\to\infty}\big[x^{n}_i(t)+M_1(t-\tau)\big]-\ve\\
		&~=~x^{*}_i(t)+M_1(t-\tau)-\ve~\leq~x_{i}^*(\tau)-{\ve\over 3},
	\end{split}
	\]
	and this yields a contradiction. Hence,   
	\[
	\lim_{n\to\infty}x_i^n(\tau)~=~x_i^*(\tau),\qquad\forall \tau\in (0,T)\backslash \mathcal{D},
	\]
	and the first equation of (\ref{inf-1}) holds for every $t\in [0,T]$. With the same argument, we  obtain the second equation of (\ref{inf-1}).
	\medskip

	{\bf  3.}  Next in order to show the absolute continuity of \((x^*,y^*)\), we claim that $x_i^*$ is locally Lipschitz continuous for all $i\in\{1,2,...,N_1\}$. If not, from (\ref{onesideL-x*})  there exist $i_0\in\{1,2,...,N_1\}$ and sequence of times $a_k<b_k$ such that  
	\[
	\lim_{k\to\infty}a_k~=~\lim_{k\to\infty}b_k~=~\bar{t}\in [0,T[\,,
	\]
	and
	\bel{eq1}
	x^*_{i_0}(b_k)-x^*_{i_0}(a_k)~\geq~k\cdot (b_k-a_k),\qquad\forall k\geq 1\,.
	\eeq
	Recalling (\ref{inf-1}), for any fixed $k\geq 1 $, there exists a  supersolution \((x^{n_k}, y^{n_k}) \) for $n_k\geq 1$ sufficiently large such that 
	\bel{eq2}
	\big|x^*_{i_0}(a_k)-x^{n_k}_{i_0}(a_k)\big|~\leq~{k\over 2}\cdot (b_k-a_k)\,.
	\eeq			
	Introduce the function $f^{\sharp}:\R\times\R^{N_1}\to \R^{N_1}$ defined by
	\bel{f-s}
	f^{\sharp}_i(t,x)~=~\max\big\{-M_1,f_i(t,x,y^{n_k}(t))\big\},\qquad (t,x)\in \R\times\R^{N_1}, i\in \{1,\cdots, N_1\}.
	\eeq
	We claim that the Cauchy problem 
	\bel{oode}
	\dot{x}(t)~=~f^{\sharp}(t,x(t)),\qquad {x}(a_k)~=~x^{n_k}(a_k), 
	\eeq
	admits a solution on $[a_k,T]$. Indeed, by Carath\'eodory's existence theorem, the problem has a local solution  $x^{\sharp,n_k}$ defined on $[a_k,T_1]$ for some $a_k<T_1<T$. For every $i\in \{1,\cdots, N_{\color{blue}1}\}$, $t\in [a_k,T_1]$, it holds 
	\[
	x^{\sharp,n_k}_i(t)~\geq~x_i^{n_k}(a_k)-M_1t~\geq~x_i^{n_k}(a_k)-M_1T.
	\]
	On the other hand, by (\ref{b-xn}) and (\ref{f-s}), we have
	\[
	\dot{x}^{n_k}(t)~\succeq~f(t,x^{n_k}(t),y^{n_k}(t))~=~f^{\sharp}(t,x^{n_k}(t)),\qquad a.e.~ t\in [0,T].
	\] 
	Using a comparison principle argument, we derive
	\bel{cp1}
	{x}^{\sharp,n_k}(t)~\preceq~x^{n_k}(t)~\preceq~(M,\cdots, M),\qquad\forall t\in [a_k,T_1],
	\eeq
	which yields 
	\[
	\begin{split}
		\max_{i\in \{1,\cdots, N_1\},t\in [a_k,T_1]} \big|x^{\sharp,n_k}_i(t)\big|&~\leq~\max\left\{M, \big|x^{n_k}_1(a_k)-M_1T\big|,\cdots, \big|x^{n_k}_{N_1}(a_k)-M_1T\big|\right\}\\
		&~\leq~M+M_1T
	\end{split}
	\]
	Hence, the solution $x^{\sharp,n_k}$ of the Cauchy problem (\ref{oode}) can be extended up to time $T$ and satisfies
	\bel{bbb}
	\max_{i\in \{1,\cdots, N_1\},t\in [a_k,T]} \big|x^{\sharp,n_k}_i(t)\big|~\leq~M+M_1T,\qquad {x}^{\sharp,n_k}(t)~\preceq~x^{n_k}(t)\qquad\forall t\in [0,T],
	\eeq
	and 
	\[
	\dot{x}^{\sharp,n_k}(t)~=~f^{\sharp}(t,x^{\sharp, n_k}(t))~\succeq~f(t,x^{\sharp, n_k}(t), y^{n_k}(t)), \qquad a.e.~t\in [0,T].
	\]
	Next,  we define  $\tilde{x}^{n_k}$ as as the absolutely continuous function
	\[
	\tilde{x}^{n_k}(t)~=~\begin{cases}
		x^{n_k}(t),\qquad\forall t\in [0,a_k),\\[2mm]
		x^{\sharp,n_k}(t),\qquad\forall t\in [a_k, T],
	\end{cases}
	\]
	which satisfies
	\[
	\dot{\tilde{x}}^{n_k}(t)~\succeq~f(t,\tilde{x}^{n_k}(t), y^{n_k}(t)), \qquad a.e.~t\in [0,T].
	\]			
	By the quasi-monotonicity  property of $g$ and (\ref{bbb}), it holds 
	\[
	\dot{y}^{n_k}(t)~\preceq~g(t,x^{n_k}(t),y^{n_k}(t))~\preceq~g(t,\tilde{x}^{n_k}(t),y^{n_k}(t)), \qquad a.e.~t\in [0,T],
	\]
	and this implies that  $(\tilde{x}^{n_k},y^{n_k})$ is also a supersolution of (\ref{1})-(\ref{2}). By (\ref{b-xn}) and (\ref{bbb}), we have  
	\[
	\sup_{t\in [a_k,T]}\big|\dot{\tilde{x}}^{n_k}(t)\big|~\leq~M_2~\doteq~M_1+\sup_{t\in [0,T], |x|\leq N_1M+N_1M_1T, |y|\leq M} |f(t,x,y)|.
	\]
	Together with (\ref{eq1}), (\ref{eq2}) and (\ref{cp1}),  by choosing $k\geq 2M_2$ , we obtain
	\[
	\begin{split}
		x^{*}_{i_0}(b_k)-\tilde{x}^{n_k}_{i_0}(b_k)&~=~\big[x^{*}_{i_0}(b_k)-x^{*}_{i_0}(a_k)\big]+\big[x^{*}_{i_0}(a_k)- x^{n_k}_{i_0}(a_k)\big]-\big[\tilde{x}^{n_k}_{i_0}(b_k)-\tilde x^{n_k}_{i_0}(a_k)\big]\\
		&~\geq~{k\over 2}\cdot (b_k-a_k)-\big[\tilde x^{n_k}_{i_0}(b_k)- \tilde x^{n_k}_{i_0}(a_k)\big]~\geq~\left({k\over 2}-M_2\right)\cdot (b_k-a_k)>0,
	\end{split}
	\]
	contradicting the minimality of $x^*$ at \(b_k\).
	
	Similarly, one can show that  $y_j^*$ is locally Lipschitz continuous for all $j\in\{1,2,...,N_2\}$.
	\medskip

	{\bf  4.} We are ready to prove that $(x^*,y^*)$ is a solution of (\ref{1})-(\ref{2}). Let $\tau\in (0,T)$ be a differentiable point of $(x^*,y^*)$. For any $i\in \{1,2,...,N_1\}$, we have
	\[
	\begin{split}
		\dot{x}_i^{*}(\tau)&~=~\lim_{s\to 0^+}~{x_i^*(\tau+s)-x_i^*(\tau)\over s}~=~\lim_{s\to 0^+}\lim_{n\to\infty}{x_i^n(\tau+s)-x^n_i(\tau)\over s}\\
		&~\geq~\ds\lim_{s\to 0^+}\lim_{n\to\infty}{1\over s}\ds\int^{\tau+s}_{\tau}f_i(t,x^n(t),y^n(t))dt~=~f_i(\tau,x^*(\tau),y^*(\tau)).
	\end{split}
	\]
	Similarly, for any $j\in\{1,2,...,N_2\}$, it holds
	\[
	y_j^*(\tau)~\leq~g_j(\tau,x^*(\tau),y^*(\tau))\,.
	\]
	To complete the proof, we show that 
	\[
	\dot{x}_i^{*}(\tau)~\leq~f_i(\tau,x^*(\tau),y^*(\tau))~~\qquad\mathrm{and}~~\qquad y_j^*(\tau)~\geq~g_j(\tau,x^*(\tau),y^*(\tau)).
	\]
	Assume that  $\dot{x}^*_{i_0}(\tau)>f_{i_0}(\tau,x^*(\tau),y^*(\tau))+\ve$ for some $\ve>0$ and $i_0\in \{1,2,...,N_1\}$. As in step 3, let $f^{*}:\R\times\R^{N_1}\to\R^{N_1}$	be such that  
	\[
	f^{*}_i(t,x)~=~\max \big\{-M_1, f_i(t,x,y^{*}(t))\big\},\qquad (t,x)\in \R\times\R^{N_1}, i\in \{1,\cdots, N_1\}.
	\]			
	We  define $x^{\tau}$ as an  absolutely continuous function given by 
	\[
	x^{\tau}(t)~=~x^*(t),\qquad\forall t\in [0,\tau]\,,
	\]
	and $x^{\tau}$ a the solution of the Cauchy problem
	\[
	\dot{x}(s)~=~f^{*}(s,x)\qquad a.e.~s\in [\tau,T].
	\]
	With the same argument in Step 3, we can show that $(x^\tau,y^*)$ is  a supersolution of (\ref{1})-(\ref{2}). Then, for $t-\tau>0$ sufficiently small we have $x_{i_0}^*(t)>x_{i_0}^{\tau}(t)$ against the minimality of $x^*$. Thus, we finally get  
	\[
	\dot{x}_i^{*}(\tau)~=~f_i(t,x^*(\tau),y^*(\tau)),\qquad a.e.~\tau\in [0,T].
	\]
	With the same argument, one can show that for any $j\in\{1,2,...,N_2\}$ it holds
	\[
	\dot{y}^*_j(\tau)~=~g_j(t,x^*(\tau),y^*(\tau)), \qquad a.e.~\tau\in [0,T],
	\]
	and  $(x^*,y^*)$ is a solution to (\ref{1})-(\ref{2}). 
	\endproof

	\begin{remark}\label{boun}
		The condition on the uniformly bounded below of  all supersolutions of (\ref{1})-(\ref{2}) in Theorem \ref{T1} is necessary. Indeed, given any constant $a,b\in \R$, consider the simple $2\times 2$ system which satisfies the quasi-monotonicity property {\bf (M)}
		\bel{ex1}
		\left\{ \bega{rl} \dot x&=~~y\\[2mm]
		\dot y&=~-x\enda\right.
		\qquad~~~\mathrm{with}\qquad~~~
		\left\{ \bega{rl}  x(0)&=~~a,\\[2mm]
		y(2\pi)&=~b.\enda\right.
		\eeq
		Assume that $t\mapsto (x(t),y(t))$ is a solution to (\ref{ex1}). Then, we have 
		\[
		{d\over dt} \big[x^2(t)+y^2(t)\big]~=~0,\qquad t\in (0,2\pi),
		\]
		and this implies 
		$$
		\left\{ \bega{rl} x(t)&=~r\,\sin(\theta+t)\,,\\[3mm]
		y(t)&=~r\,\cos(\theta+t)\,.\enda\right.
		\qquad
		$$
		By the boundary conditions, we get
		\[
		r~=~\sqrt{a^2+b^2},\qquad\qquad \theta~=~\arcsin\left({a\over \sqrt{a^2+b^2}}\right).
		\]

		Now choose any  $a^*\geq a$ and $b^*\geq b$. Let  $(x^*, y^*)$ be
		the corresponding solution to  the boundary value problem for (\ref{ex1}) with $(a,b)=(a^*,b^*)$. 
		By definition this is a supersolution of the original problem (\ref{ex1}).    Observe that
		\[
		\min_{t\in [0,\pi]}x^*(t)~=~-\sqrt{[a^*]^2+[b^*]^2},
		\]
		the family of supersolutions of  the boundary value problem (\ref{ex1}) is not uniformly bounded below. Hence, Theorem \ref{T1} can not be applied here. 
		Moreover,  it is  not true that
		$$(x(t),y(t))~\preceq~(x^*(t),y^*(t))\qquad\forall t\in [0, 2\pi].$$
		Hence, the unique solution $(x,y)$ of (\ref{ex1}) is not a minimal solution of (\ref{ex1}). Moreover, in this case,  a standard comparison principle does not hold.
	\end{remark}

	Next, we shall provide additional assumptions on $f$ and $g$ to ensure that  the family of supersolutions of  the boundary value problem (\ref{1})-(\ref{2}) is uniformly bounded below. In order to do so,  let's introduce ${\bf f}_{\min},{\bf g}_{\max}:\R^3\to\R$ which are defined as follows: for every $s,\tau\in \R$, set $x_s=(s,\cdots, s)\in \R^{N_1}$ and $y_\tau=(\tau,\cdots, \tau)\in \R^{N_2}$, it holds
	\bel{FG-1}
	\begin{cases}
		{\bf f}_{\min}(t,s,\tau)~=~\min\{f_1(t,x_s,y_\tau),\cdots, f_{N_1}(t,x_s,y_\tau) \},\\[2mm]
		{\bf g}_{\max}(t,s,\tau)~=~\max\{g_1(t,x_s,y_\tau),\cdots, g_{N_2}(t,x_s,y_\tau)\}.
	\end{cases}
	\eeq

	It is clear that both ${\bf f_{\min}, g_{\max}}:\R^3\to\R$ are continuous w.r.t $s,\tau$ and  measurable w.r.t $t$. Moreover, the system of ODEs
	\bel{ODE-XY}
	\begin{cases}
		\dot{{\bf x}}~=~{\bf f}_{\min}(t,{\bf x},{\bf y}),\\[2mm]
		\dot{\bf y}~=~{\bf g}_{\max}(t,{\bf x},{\bf y}),
	\end{cases}
	\quad\qquad a.e.~t\in [0,T],
	\eeq
	satisfies the  quasi-monotonicity property {\bf (M)} where ${\bf x, y}:\R \to \R$. 
	\begin{lemma}\label{Min-XY}
		For every  super-solution $(x(t),y(t))$  of  (\ref{1})-(\ref{2}), we define 
		\bel{XY}
		{\bf x}(t)~=~\min_{i\in\{1,\cdots,{N_1}\}} x_i(t),\qquad {\bf y}(t)~=~\min_{j\in \{1,\cdots, {N_2}\}}y_j(t),\qquad t\in [0,T].
		\eeq
		Then the map $t\mapsto ({\bf x}(t),{\bf y}(t))$ is a super-solution of (\ref{ODE-XY})  with boundary data
		\bel{Ini-XY}
		{\bf x}(0)~=~\min\{\bar{x}_1,\cdots,\bar{x}_{N_1}\},\qquad {\bf y}(T)~=~\min\{\bar{y}_1,\cdots,\bar{y}_{N_2}\}.
		\eeq
	\end{lemma}
	{\bf Proof.} Since $x_i$'s and $y_j$'s are absolutely continuous for all $i\in \{1,\cdots, {N_1}\}$ and $j\in \{1,\cdots, {N_2}\}$, we have that  both ${\bf x}$ and ${\bf y}$ are also absolutely continuous. For every $t\in (0,T)$ a differentiable point of $x_i,{\bf x},y_j,{\bf y}$, assume that ${\bf x}(t)=x_{i_0}(t)$ and ${\bf y}(t)=y_{j_0}(t)$ for some $i_0\in\{1,\cdots, {N_1}\}$, $j_0\in\{1,\cdots, {N_2}\}$. Since  $(x(t),y(t))$  is a  super-solution of  (\ref{1})-(\ref{2}), we derive
	\bel{Xe}
	\begin{split}
		\dot{{\bf x}}(t)&~=~\lim_{s\to 0-}{{\bf x}(t+s)-{\bf x}(t)\over s}~\geq~\lim_{s\to 0-} {x_{i_0}(t+s)-x_{i_0}(t)\over s}~=~\dot{x}_{i_0}(t)~\geq~f_{i_0}(t,x(t),y(t)),
	\end{split}
	\eeq
	\bel{Ye}
	\dot{{\bf y}}(t)~=~\lim_{s\to 0+}{{\bf y}(t+s)-{\bf y}(t)\over s}~\leq~\lim_{s\to 0+} {y_{j_0}(t+s)-y_{j_0}(t)\over s}~=~\dot{y}_{j_0}(t)~\leq~g_{j_0}(t,x(t),y(t)).
	\eeq
	By (\ref{FG-1}), (\ref{XY}), and  the quasi-monotonicity property {\bf (M)} of (\ref{1}), we derive 
	\[
	\begin{split}
		f_{i_0}(t,x(t),y(t))&~=~f_{i_0}(t,x_1(t),\cdots, x_{i_0-1}(t), {\bf x}(t), x_{i_0+1}(t),\cdots x_{N_1}(t),y_1(t),\cdots, y_{N_2}(t))\\
		&~\geq~f_{i_0}(t, {\bf x}(t), \cdots,{\bf x}(t),{\bf y}(t),\cdots,{\bf y}(t))\\
		&~\geq~{\bf f}_{\min}(t,{\bf x}(t),{\bf y}(t)),
	\end{split}
	\]
	and 
	\[
	\begin{split}
		g_{j_0}(t,x(t),y(t))&~=~g_{j_0}(t,x_1(t),\cdots, x_{N_1}(t),y_1(t),\cdots, y_{i_0-1}(t), {\bf y}(t), y_{i_0+1}(t),\cdots, y_n(t))\\
		&~\leq~g_{j_0}(t, {\bf x}(t), \cdots,{\bf x}(t),{\bf y}(t),\cdots,{\bf y}(t))\\
		&~\leq~{\bf g}_{\max}(t,{\bf x}(t),{\bf y}(t)).
	\end{split}
	\]
	Hence, from (\ref{Xe})-(\ref{Ye}), we get 
	\[
	\dot{{\bf x}}(t)~\geq~{\bf f}_{\min}(t,{\bf x}(t),{\bf y}(t)),\qquad \dot{{\bf y}}(t)~\leq~{\bf g}_{\max}(t,{\bf x}(t),{\bf y}(t)),
	\]
	this implies that the map $t\mapsto ({\bf x}(t),{\bf y}(t))$ provides a super-solution of (\ref{ODE-XY}),(\ref{Ini-XY}).
	\endproof

	As a consequence, we obtain the following corollary. 
	
	\begin{corollary} Assume that the family of supersolutions of  the system of ODEs   (\ref{ODE-XY}) with boundary condition (\ref{Ini-XY}) is uniformly bounded below. Then the family of supersolutions of  the boundary value problem of  (\ref{1})-(\ref{2}) is also uniformly bounded below.
	\end{corollary}

	Now, we  establish sufficient condition on ${\bf f}_{\min},{\bf g}_{\max}$ which ensure that the family of supersolutions of  the system of ODEs   (\ref{ODE-XY}) with boundary condition (\ref{Ini-XY}) is  uniformly bounded below. To obtain that this family is nonempty, we also need to impose  another conditions on   functions ${\bf f}_{\max},{\bf g}_{\min}:\R^3\to\R$ which are defined as follows: for every $s,\tau\in \R$, set $x_s=(s,\cdots, s)\in \R^{N_1}$ and $y_\tau=(\tau,\cdots, \tau)\in \R^{N_2}$, it holds
	\bel{FG}
	\begin{cases}
		{\bf f}_{\max}(t,s,\tau)~=~\max\{f_1(t,x_s,y_\tau),\cdots, f_{N_1}(t,x_s,y_\tau) \},\\[2mm]
		{\bf g}_{\min}(t,s,\tau)~=~\min\{g_1(t,x_s,y_\tau),\cdots, g_{N_2}(t,x_s,y_\tau)\}.
	\end{cases}
	\eeq
	
	The following theorem establishes the existence of a minimal solution to \eqref{1}–\eqref{2}, provided certain assumptions are satisfied for \({\bf f}_{\min},{\bf f}_{\max},{\bf g}_{\min},{\bf g}_{\max}.\)

	\begin{theorem}\label{Main2} Under the  quasi-monotonicity property {\bf (M)} and {\bf (A)},  assume that one of the following conditions holds:
		\begin{itemize}
			\item [(i).] There exist two functions $\alpha_1,\alpha_2:[0,T]\times \R\to\R$ such that for all  $(t,s,\tau)\in [0,T]\times\R^2$
			\bel{as1}
			\alpha_1(t,s)~\leq~{\bf f}_{\min}(t,s,\tau)~\leq~{\bf f}_{\max}(t,s,\tau)~\leq~\alpha_2(t,s).
			\eeq
			For every $a,b, M\in\R$, the Cauchy problems 
			\bel{cau1}
			\dot{\gamma}(t)~=~\alpha_{i}(t,\gamma),\qquad \gamma(0)~=~a,\qquad i=1,2,
			\eeq
			\[
			\dot{\eta}(t)~=~{\bf g}_{\min}(t,M,\eta),\qquad \eta(T)~=~b,
			\]
			admit solutions defined on  $[0,T]$.
			\item [(ii).] There exists two functions $\beta_1,\beta_2:[0,T]\times \R\to\R$ such that $(t,s,\tau)\in [0,T]\times\R^2$,
			\[
			\beta_1(t,\tau)~\leq~{\bf g}_{\min}(t,s,\tau)~\leq~{\bf g}_{\max}(t,s,\tau)~\leq~\beta_2(t,\tau).
			\]
			For every $a,b, M\in\R$, the Cauchy problems 
			\[
			\dot{\gamma}(t)~=~\beta_i(t,\gamma),\qquad \gamma(T)~=~a,\qquad i=1,2,
			\]
			\[
			\dot{\eta}(t)~=~{\bf f}_{\max}(t,\eta,M),\qquad \eta(0)~=~b,
			\]
			admit solutions defined on  $[0,T]$.
		\end{itemize}
		Then  the boundary value problem  (\ref{1})-(\ref{2}) admits a minimal solution for every $(x,y)\in \R^{N_1}\times\R^{N_2}$.
	\end{theorem}
	{\bf Proof.} We will only show that our assertion holds under condition (i), as the proof for condition (ii) is similar.

	{\bf 1.} For every  super-solution  $t\mapsto ({\bf x}(t),{\bf y}(t))$  of    the two point boundary problem (\ref{ODE-XY})  with boundary conditions $({\bf x}(0),{\bf y}(T))=(\bar{\bf x},\bar{\bf y})$, it holds
	\[
	\dot{\bf x}(t)~\geq~{\bf f}_{\min}(t,{\bf x},{\bf y})~\geq~\alpha_1(t,{\bf x}(t)),\qquad\quad \bar{\bf x}(0)~\geq~\bar{\bf x}.
	\]
	Let $\gamma_1(\cdot)$ be a solution of the Cauchy problem 
	\[
	\dot{\gamma}(t)~=~\alpha_{1}(t,\gamma),\qquad \gamma(0)~=~\bar{\bf x}.
	\]
	By a standard comparison principle argument, we obtain that 
	\[
	{\bf x}(t)~\geq~\gamma_1(t)\qquad\forall  t\in [0,T].
	\]
	Set $\ds\gamma_{\min}=\min_{t\in [0,T]}\gamma_1(t)$. The quasi-monotonicity property of (\ref{ODE-XY}) yields 
	\[
	\dot{\bf y}(t)~\leq~{\bf g}_{\max}(t,{\bf x},{\bf y})~\leq~{\bf g}_{\max}(t,\gamma_{\min},{\bf y})
	,\qquad  {\bf y}(T)~\geq~\bar{\bf y}.
	\]
	Again, by a standard comparison principle argument and the assumption (i), ${\bf y}(t)$ is bounded below by a solution $\eta_1$ of the Cauchy problem
	\[
	\dot{\eta}(t)~=~{\bf g}_{\max}(t,\gamma_{\min},\eta),\qquad \eta(T)~=~\bar{\bf y}. 
	\]
	Hence, every  super-solution  $t\mapsto ({\bf x}(t),{\bf y}(t))$  of    the two point boundary problem (\ref{ODE-XY})  with boundary conditions $({\bf x}(0),{\bf y}(T))=(\bar{\bf x},\bar{\bf y})$ is bounded below by $\min\left\{\gamma_{\min},\min_{t\in [0,T]}\eta_1(t)\right\}$. By Lemma \ref{Min-XY},  the family of supersolutions of  the boundary value problem (\ref{1})-(\ref{2}) is uniformly bounded from below. 
	\medskip

	\n {\bf 2. } Next, for every given boundary data $(\bar{x},\bar{y})\in \R^{N_1}\times \R^{N_2}$, we set 
	\bel{Ib}
	\bar{\bf x}_{\max}~=~\max\{\bar{x}_1,\cdots, \bar{x}_{N_1}\},\qquad \bar{{\bf y}}_{\max}~=~\max\{\bar{y}_1,\cdots, \bar{y}_{N_2}\}.
	\eeq
	Let $t\mapsto \gamma_2(t)$ be solution of the Cauchy problem 
	\[
	\dot{\gamma}(t)~=~\alpha_{2}(t,\gamma),\qquad \gamma(0)~=~\bar{\bf x}_{\max}.
	\]
	Set $\ds\gamma_{\max}=\max_{t\in [0,T]}\gamma_2(t)$. Let $t\mapsto \eta_2$ be the solution of the  Cauchy problem 
	\[
	\dot{\eta}(t)~=~{\bf g}_{\max}(t,\gamma_{\max},\eta),\qquad \eta(T)~=~ \bar{{\bf y}}_{\max}.
	\]
	We claim that the map  $t\mapsto (\tilde{x}(t), \tilde{y}(t))\in \R^{N_1}\times\R^{N_2}$ defined by 
	\[
	\tilde{x}_i(t)~=~\gamma_2(t),\qquad \tilde{y}_j(t)~=~\eta_1(t),\qquad i\in \{1,\cdots, {N_1}\}, j\in \{1,\cdots, {N_2}\}
	\]
	provide a supersolution of (\ref{1})-(\ref{2}). Indeed, by (\ref{Ib}), we have 
	\[
	\bar{x}~\preceq~\tilde{x}(0)\qquad~~\mathrm{and}~~\qquad\bar{y}~\preceq~\tilde{y}(T).
	\]
	For every $i\in \{1,\cdots, {N_1}\}$ and  $j\in \{1,\cdots, {N_2}\}$, by (\ref{as1}) and  the  quasi-monotonicity property {\bf (M)} of (\ref{1}), we derive
	\[
	\begin{split}
		\dot{\tilde{x}}_i(t)&~=~\dot{\gamma}_2(t)~=~\alpha_2(t,\gamma_2(t))~\geq~{\bf f}_{\max}(t,\gamma_2(t),\eta_2(t))~\geq~f_i(t,\tilde{x}(t),\tilde{y}(t)),
	\end{split}
	\]
	\[
	\dot{\tilde{y}}_i(t)~=~{\bf g}_{\min}(t,\gamma_{\max},\eta_2(t))~\leq~g_j(t,\gamma_{\max},\cdots, \gamma_{\max}, \tilde{y}(t))~\leq~g_j(t,\tilde{x}(t),\tilde{y}(t)),
	\]
	and this yields 
	\[
	\dot{\tilde{x}}(t)~\succeq~ f(t,\tilde{x},\tilde{y}),\qquad\qquad \dot{\tilde{y}}(t)~\preceq~ g(t,\tilde{x},\tilde{y}).
	\]
	Since the family of supersolutions of  the boundary value problem (\ref{1})-(\ref{2}) is nonempty and  uniformly bounded below, by Theorem \ref{T1} we conclude that (\ref{1})-(\ref{2}) admits a minimal solution for every $(x,y)\in \R^{N_1}\times\R^{N_2}$.
	\endproof
	\v

	We conclude this section with an application of Theorem \ref{Main2} to a two-point boundary value problem for a Hamiltonian system.

	From now on, by $\C^k(\R^d)$ we denote the Banach space of all bounded functions with bounded, continuous  partial derivatives up to order $k$, see for example  \cite{BLN, Evans}.
	\begin{remark}\label{r1}
		Note that if both $f$ and $g$ are $\mathcal{C}^1$ smooth w.r.t. $x,y$ and  satisfy
		\bel{QM}\bega{rl}\ds
		{\partial f_i\over\partial x_k}~\geq~ 0\,,\qquad\forall k\not= i,
		\qquad\qquad &\ds {\partial f_i\over\partial y_k}~\geq~ 0,\qquad\forall k=1,\ldots,{N_1},\\[4mm]\ds
		{\partial g_j\over\partial x_k}~\leq~ 0,\qquad\forall k\,,
		\qquad\qquad  &\ds{\partial g_j\over\partial y_k}~\leq~ 0\,,
		\qquad\forall k\not= j,\enda\eeq
		then the system of ODEs (\ref{1}) has  the  quasi-monotonicity property {\bf (M)}. 
	\end{remark}

	In general, the two-point boundary value Hamiltonian system arising from optimal control problems may fail to admit a solution (see, e.g., \cite[Example 2.1]{BN1}). Nevertheless, combining Remark \ref{r1} with Theorem \ref{Main2} yields the following corollary:

	\begin{corollary}
		[Hamiltonian system] Consider the Hamiltonian system
		\bel{HJ}
		\begin{cases}
			\dot{x}(t)~=~D_pH(x,p),\\[2mm]
			\dot{p}(t)~=~-D_xH(x,p),
		\end{cases}
		\qquad\qquad
		\begin{cases}
			x(0)~=~\bar{x},\\[2mm]
			p(T)~=~\bar{p},
		\end{cases}
		\eeq
		with the Hamiltonian $H\in\mathcal{C}^2(\R^d\times\R^d,\R)$ satisfying 
		\bel{AH}
		\partial_{x_ix_j}H(x,p), \partial_{p_ip_j}H(x,p), \partial_{x_ip_k}H(x,p)~\geq~0
		\eeq
		for all $i,j,k\in \{1,\cdots, d\}$ with $i\neq k$. The Hamiltonian system (\ref{HJ}) admits a minimal solution.
	\end{corollary}
	{\bf Proof.} By Remark \ref{r1} and (\ref{AH}),  the system of ODEs  (\ref{HJ}) has the  quasi-monotonicity property {\bf (M)}. Moreover, since $H$ is in $\mathcal{C}^2(\R^d\times\R^d,\R)$, both conditions (i) and (ii) in Theorem \ref{Main2} are satisfied. Therefore, by Theorem \ref{Main2}, the Hamiltonian system (\ref{HJ}) admits a minimal solution.
	\endproof
	
	\section{An application to mean field games}
	
	In this section, we use  Theorem \ref{T1} to establish  a non-uniqueness  result on strong solutions to a class of   mean field games  with a continuum of players. More precisely, let   $V:\R^d\to \R$ be a  $\mathcal{C}^2$ potential function such that  
	\bel{con-V}
	DV(0)~=~0, \qquad V_{x_ix_j}(x)~\geq~0\quad\forall  i,j\in \{1,\cdots, d\}, x\in \R^d.
	\eeq
	We consider a mean field game where each player $\xi\in [0,1]$ minimizes the same cost
	\bel{Cost}
	J[u]~=~\int_{0}^{T}{|u(t)|^2\over 2}+V(x(t))+\kappa \cdot |x(t)-b(t)|^2dt,
	\eeq
	subject to  a controlled ODE
	\bel{dy1}
	\dot{x}~=~u~\in~\R^d,\qquad x(0,\xi)~=~0,\qquad\forall \xi\in [0,1].
	\eeq
	Here, $b$ denotes the barycenter of the distribution of players
	\bel{ba}
	b(t)~=~\int_{0}^{1}x(t,\xi)d\xi.
	\eeq

	\begin{definition}[Strong solution to MFG]
		We say that a family of control functions $t\mapsto u(t,\xi)\in\R^d$ and corresponding trajectories $t\mapsto x(t,\xi)\in\R^d$ defined for all $\xi\in [0,1]$, $t\in [0,T]$ is a strong solution to the mean field game (\ref{Cost})-(\ref{ba}) if the following holds:
		\medskip
		
		For a.e. $\xi\in [0,1]$ the control $u(\cdot, \xi)$ and the trajectory $x(\cdot,\xi)$ provide an optimal solution to  the optimal control problem (\ref{Cost})-(\ref{dy1}) with $b$ defined in (\ref{ba}).
	\end{definition}

	By the above definition, a mean field game thus yields a (possibly multivalued) $b\mapsto \Phi(b)$ from $\mathcal{C}([0,T],\R^d)$ into itself. Namely, given $b\in \mathcal{C}([0,T],\R^d) $, for each player $\xi\in [0,1]$  consider an optimal trajectory $t\mapsto x^b(t,\xi)$  of the corresponding optimal control problem  (\ref{Cost})-(\ref{dy1}). Then we set 
	\[
	\Phi(b)(t)~=~\int_{0}^{1}x^b(t,\xi)d\xi,\qquad t\in [0,T].
	\]
	By definition, a fixed point of this composed map  
	\bel{mot}\bega{c}\qquad\qquad~~~~ b(\cdot)~~\quad \mapsto~~\quad \bigl\{ x^b(\cdot, \xi)\,;~\xi\in \Omega\bigr\}
	\qquad\mapsto\qquad \Phi(b)~\doteq~\ds\int_{0}^{1}x^b(t,\xi)d\xi\\[3mm]
	[\hbox{moments}]\quad~~\mapsto~~\quad [\hbox{optimal trajectories}]\quad~\mapsto~\quad [\hbox{moments}]
	\enda
	\eeq
	yields a strong solution to the mean field game (\ref{Cost})--(\ref{ba}).
	\medskip

	In general, the map $\Phi$ can  be multivalued. Indeed, for some $b$, there can be a subset $I\subset [0,1]$  with positive measure, such that each player in $\xi\in I$  has two or more optimal
	trajectories. For this reason, a mean field game may not have a solution in the strong sense
	considered in the above definition.

	Following the standard literature on fixed points of continuous or multivalued maps, we introduce the concept of stable  solutions of a mean field game.
	
	\begin{definition} Let $x^*=x^*(t,\xi)$ be  a strong solution to the mean field game (\ref{Cost})-(\ref{ba}) such that 
		\[
		b^*(t)~=~\int_{0}^{t}x^*(t,\xi)d\xi,\qquad t\in [0,T],
		\]
		is a fixed point of $\Phi$. We say that the solution $x^*$ is stable if the corresponding function  $b^*$ is  a stable fixed point of the multivalued map $b\mapsto \Phi(b)$ in  $\mathcal{C}([0,T],\R^d)$. Namely, for every $\ve>0$, there is $\delta>0$ such that for every sequence $b_n$ such that 
		\[
		\|b_0-b^*\|_{\mathcal{C}^0}~\leq~\delta,\qquad b_{n}~\in~\Phi(b_{n-1}),\qquad n\geq 1,
		\]
		it holds that  $\|b_n-b^*\|_{\mathcal{C}^0}\leq \ve$ for all $n\geq 1$. In addition, if every such sequence $b_n$ converges to $b^*$, then we say that the solution is asymptotically stable.
		
	\end{definition}

	Using Theorem \ref{Main2} we show that the mean field game   (\ref{Cost})-(\ref{ba}) admits at least two stable solutions. 
	\begin{theorem}\label{Main3} Assume that $V\in\mathcal{C}^2(\R^d)$ satisfies (\ref{con-V}) and 
		\bel{conV2}
		0~<~\gamma_{\min}~\doteq~\min_{i\in \{1,\cdots, d\}}\left(\sum_{j=1}^{d}V_{x_ix_j}(0)\right)~\leq~\gamma_{\max}~\doteq~\max_{i\in \{1,\cdots, d\}}\left(\sum_{j=1}^{d}V_{x_ix_j}(0)\right).
		\eeq
		Then  the following hold:
		\begin{itemize}
			\item [(i).] For every $\kappa,T>0$ such that 
			\bel{conT}
			\kappa~\geq~\big\|D^2V\big\|_{\infty}, \qquad T^{2/3}~\geq~\max\left\{{\gamma^2_{\max}\over 8 \gamma_{\min}},\left({1\over \gamma_{\min}}+\sqrt{1+{1\over \gamma^2_{\min}}}\right)^{3/2}\right\},
			\eeq
			the mean field game (\ref{Cost})-(\ref{ba}) admits at least two solutions. These have
			the form
			\[
			x^j(t,\xi)~=~y_j(t),\qquad\qquad j=1,2,
			\]
			with $y_1\equiv 0$, while $y_2\neq 0$ is monotone increasing.
			\item [(ii).] The zero  solution $x^1$ is asymptotically stable and the negative solution $x^2$ is stable. 
		\end{itemize}
	\end{theorem}
	
	{\bf Proof.} {\bf 1.} For a given  barycenter $b(\cdot)$, the Pontryagin necessary condition of (\ref{Cost})-(\ref{dy1}) takes 
	the form
	\bel{ex2}
	\left\{ \bega{rl} \dot x(t)&=~~-p(t)\\[2mm]
	\dot p(t)&=~\ds-DV(x(t))-2\kappa (x(t)-b(t))\enda\right.
	\eeq
	with boundary conditions
	\bel{bd2}
	x(0)~=~0,\qquad\qquad p(T)~=~0.
	\eeq
	Consider the Hamiltonian associated to (\ref{Cost})-(\ref{ba})
	\[
	H^b(x,u,p,t)~=~{|u|^2\over 2}+V(x)+\kappa\cdot |x-b(t)|^2+p\cdot u.
	\]
	The reduced Hamiltonian is given by 
	\[
	\begin{split}
		\Hat{H}^{b}(x,p,t)&~=~\min_{u\in \R} H^b(x,u,p,t)~=~V(x)+\kappa\cdot |x-b(t)|^2-{p^2\over 2}.
	\end{split}
	\]
	We compute that 
	\[
	D_x^2\Hat{H}^{b}(x,p,t)~=~\kappa\cdot \mathbb{I}_d +D^2V(x).
	\]
	Hence, for $\kappa> \big\|D^2V\big\|_{\infty}$, we have 
	\[
	D_x^2\Hat{H}^{b}(x,p,t)~>~0\qquad\forall (x,p,t)\in\R^d\times\R^d\times[0,T],
	\]
	and the map $x\mapsto \Hat{H}^{b}(x,p,t)$ is strictly convex. This  implies that the optimal control problem (\ref{Cost})-(\ref{ba}) admits a unique pair $\big(u^{b},x^{b})$ of optimal control and trajectory.  Indeed, let $(u_1,x_1)$ be another pair of optimal control and optimal trajectory  of (\ref{Cost})-(\ref{dy1}). By the Pontryagin maximum principle, there exist $p^b,p_1\in \mathcal{C}^0([0,T])$ such that both $\big(u^{b},p^b,x^{b})$ and $\big(u_{1},p_1,x_{1})$ solve (\ref{ex2}) and 
	\[
	u_1(t)~=~p_1(t),\qquad u^b(t)~=~p^b(t)\qquad a.e.~ t\in [0,T].
	\]
	Assume that $(u_1,x_1)$ and $(x^b,u^b)$ are different. Then by the strict convexity of   $x\mapsto \Hat{H}^{b}(x,p,t)$ and the optimality condition 
	we estimate
	\[
	\begin{split}
		0&~=~ \int_{0}^{T}{|u_1|^2\over 2}+V(x_1)+\kappa \cdot |x_1-b|^2-\left({|u^b|^2\over 2}+V(x^b)+\kappa \cdot |x^b-b|^2\right)dt\\
		&~=~\int_{0}^{T}\Hat{H}^{b}(x_1(t),u_1(t),p^{b}(t),t)-\Hat{H}^{b}(x^b(t),u^b(t),p^{b}(t),t)-p^b(t)\cdot (u_1(t)-u^b(t))dt \\
		&~\geq~ \int_{0}^{T}\Hat{H}^{b}(x_1(t),p^{b}(t),t)-\Hat{H}^{b}(x^b(t),p^{b}(t),t)-p^b(t)\cdot \big(\dot{x}_1-\dot{x}^b(t)\big)dt\\
		&~>~\int_{0}^{T}\partial_x\Hat{H}^{b}(x^b(t),p^{b}(t),t)\cdot \big(x_1(t)-x^b(t)\big)-p^b(t)\cdot \big(\dot{x}_1-\dot{x}^b(t)\big)dt\\
		&~=~\int_{0}^{T}-\dot{p}^b(t)\cdot \big(x_1(t)-x^b(t)\big)-p^b(t)\cdot \big(\dot{x}_1-\dot{x}^b(t)\big)dt=\int_{0}^{T}{d\over dt}\big[p^b(t)\cdot \big(x_1(t)-x^b(t)\big)\big]\\
		&~=~p^b(0)\cdot \big(x_1(0)-x^b(0)\big)-p^b(T)\cdot \big(x_1(T)-x^b(T)\big)~=~0,
	\end{split}
	\]
	and this yields a contradiction. Hence, all the optimal trajectories $x(\cdot,\xi)$ of the mean field game coincide for $\xi\in [0,1]$. From (\ref{ex2}), the trajectories 
	$x(\cdot,\xi)=b(\cdot)$  provided a strong solution to the mean field game (\ref{Cost})-(\ref{ba}) if and only if  there exists $p(t)$ such that $(b,p)$ is a solution of the two point boundary problem
	\bel{ex3}
	\left\{ \bega{rl} \dot x&=~~-p,\\[2mm]
	\dot p&=~\ds-DV(x)\enda\right.
	\qquad~~~\mathrm{with}\qquad~~~
	\left\{ \bega{rl}  x(0)&=~0,\\[3mm]
	p(T)&=~0.\enda\right.
	\eeq
	
	Since $DV(0)=0$, we have that $(x,p)\equiv (0,0)$ is a trivial solution to (\ref{ex3}). This provides the first solution of the mean field game
	\[
	x^1(t,\xi)~=~y_1(t)~\doteq~0\qquad\forall \xi \in [0,1],~t\in [0,T].
	\]

	\n {\bf 2.} We claim that the two point boundary problem  (\ref{ex3})  admits another nontrivial solution. Indeed,  by the assumption (\ref{con-V}) and the smoothness of $V$, the system of ODEs in (\ref{ex3}) has  the  quasi-monotonicity property {\bf (M)} and satisfies conditions in Theorem  \ref{Main2}. Hence, the two point boundary problem (\ref{ex3}) admits a minimal solution $(\bar{x}^*(t),\bar{p}^*(t))$ and this provide the second solution of the mean field game
	\[
	x^2(t,\xi)~=~\bar{x}^*(t)~\doteq~y_2(t)\qquad\forall \xi \in [0,1],~t\in [0,T].
	\]
	To verify that $x^2(t,\cdot)$ is a nontrivial solution, we shall construct a negative super-solution of (\ref{ex3}). From the assumption (\ref{conV2}), it holds
	\bel{e}
	\gamma_{\min}\cdot (1,\cdots, 1)^{\dagger}~\preceq~{\bf e}~\doteq~D^2V(0)(1,\cdots, 1)^{\dagger}~\preceq~\gamma_{\max}\cdot (1,\cdots, 1)^{\dagger}.
	\eeq
	
	For every $\theta,\lambda \in (0,1)$ and ${\bf h}\in \R^d$, we consider the absolutely continuous map $t\mapsto \big(x^{\theta}(t),p^{\theta}(t)\big)\in \R^d\times\R^d$ defined by 
	\bel{x}
	x^{\theta,\lambda,{\bf h}}(t)~=~\begin{cases}
		-\theta t {\bf e},\qquad t\in [0,\lambda T],\\[2mm]
		-\theta \lambda T{\bf e}+\theta{\bf h}(t-\lambda T),\qquad t\in [\lambda T,T],
	\end{cases}
	\eeq
	and 
	\bel{p}
	p^{\theta, \lambda,{\bf h}}(t)~=~\begin{cases}
		-\theta {\bf e},\qquad t\in [0,\lambda T],\\[2mm]
		-\theta {\bf e}-\ds\int_{\lambda T}^{t}DV(-\theta \lambda T{\bf e}+\theta {\bf h}(\tau-\lambda T)))d\tau,\qquad t\in [\lambda T,T],
	\end{cases}
	\eeq
	such that 
	\bel{bd}
	x^{\theta, \lambda,{\bf h}}(0)~=~0,\qquad p^{\theta, \lambda,{\bf h}}(T)~=~-\theta {\bf e}-\ds\int_{\lambda T}^{T}DV(-\theta \lambda T{\bf e}+\theta {\bf h}(\tau-\lambda T)))d\tau,
	\eeq
	and 
	\bel{dv}
	\big(\dot{x}^{\theta, \lambda, {\bf h}}(t),\dot{p}^{\theta, \lambda,{\bf h}}(t)\big)~=~\begin{cases}\big(p^{\theta,\lambda, {\bf h}}(t),0\big),& t\in (0,\lambda T), \\[2mm]
		\big(\theta {\bf h},-DV(x^{\theta,\lambda, {\bf h}}(t))\big),& t\in (\lambda T,T).
	\end{cases}
	\eeq
	By (\ref{con-V}) and (\ref{bd}), for every $t\in (0,\lambda T)$ we compute 
	\[
	\begin{split}
		-DV(x^{\theta,\lambda,  {\bf h}}(t))&~=~-DV(- \theta t {\bf e})+DV(0)~=~ \int_{0}^{1}D^2V(-\tau \theta t {\bf e})(\theta t{\bf e})d\tau~\succeq~0~=~ \dot{p}^{\theta,\lambda,{\bf h}}(t),
	\end{split}
	\]
	which implies  that the map $t \mapsto \big(x^{\theta,\lambda, {\bf h}}(t),p^{\theta,\lambda, {\bf h}}(t)\big)$ satisfies the conditions for a super-solution of (\ref{ex3}) in $[0,\lambda T]$.  
	\medskip

	{\bf 3.} Next we will show that the map $t \mapsto \big(x^{ \theta, \lambda, {\bf h} }(t),p^{\theta, \lambda, {\bf h} }(t)\big)$ satisfies the condition of a super-solution of (\ref{ex3}) in $[\lambda T,T]$ for some  ${\bf h}\succeq 0$ and  $\theta>0$ sufficiently small.  By (\ref{dv}), we only need to find $\theta>0$ and ${\bf h}\succeq 0$ such that  
	\bel{dk}
	0~\preceq~p^{\bf \theta, \lambda, {\bf h}}(T)~\preceq~\sup_{t\in [\lambda T,T]}p^{\theta, \lambda, {\bf h}}(t)~\preceq~\theta{\bf h}.
	\eeq
	Recalling (\ref{p}), we write   
	\[
	\begin{split}
		p^{\theta, \lambda, {\bf h}}(t)&~=~-\theta{\bf e}-\int_{\lambda T}^{t}\big[DV(-\theta \lambda T{\bf e}+\theta{\bf h}(\tau-\lambda T)))-DV(0)\big]~d\tau\\
		&~=~-\theta{\bf e}-\int_{\lambda T}^{t}D^2V(0)(-\theta \lambda T{\bf e}+\theta {\bf h}(\tau-\lambda T)))d\tau+ E(t,\theta,{\bf h})\\
		&~=~\theta\cdot\left[-{\bf e}+\lambda T(t-\lambda T)\cdot D^2V(0)( {\bf e})-{(t-\lambda T)^2\over 2}D^2V(0)({\bf h})\right]+E(t,\theta,{\bf h}),
	\end{split}
	\]
	where the small term $E(t,\theta,{\bf h})$ satisfies
	\[
	\sup_{t\in [\lambda T,T]}|E(t,\theta,{\bf h})|~\leq~\O(1)\cdot \theta^2.
	\]
	
	From (\ref{e}), we choose 
	\bel{h}
	{\bf h}~=~{2\lambda {\bf e}\over (1-\lambda)} - {2(1+\sqrt{\theta})\over (1-\lambda)^2 T^2}\cdot (1,\cdots, 1)^{\dagger}.
	\eeq
	such that 
	\bel{cho-h}
	D^2V(0)\big(\lambda {\bf e}-(1-\lambda){\bf h}/2\big)~=~{1+\sqrt{\theta}\over (1-\lambda)T^2}D^2V(0)(1,\cdots, 1)^{\dagger}~=~{1+\sqrt{\theta}\over (1-\lambda)T^2} {\bf e}.
	\eeq
	We compute
	\[
	\begin{split}
		p^{\theta,\lambda,{\bf h}}(t)&~=~\theta\cdot\left[-{\bf e}+\lambda T(t-\lambda T)\cdot D^2V(0)( {\bf e})-{(t-\lambda T)^2\over 2}D^2V(0)({\bf h})\right]+E(t,\theta,{\bf h})\\
		&~=~\theta\cdot\left[\left({(1+\sqrt{\theta})(t-\lambda T)^2\over (1-\lambda)^2T^2}-1\right)\cdot {\bf e}+{\lambda(t-\lambda T) (T-t)\over (1-\lambda)}D^2V(0)({\bf e})\right]+E(t,\theta,{\bf h})
	\end{split}
	\]
	In particular, by (\ref{cho-h}) we have 
	\[
	p^{\theta,\lambda, {\bf h}}(T)~=~\theta\cdot\Big[-{\bf e}+(1+\sqrt{\theta})\cdot {\bf e}\Big]+E(T,\theta,{\bf h})~=~\theta^{3/2}{\bf e}+E(T,\theta,{\bf h})~\succeq~0.
	\]
	and this yields 
	\[
	p^{\theta,\lambda, {\bf h}}(T)~\succeq~0\qquad\mathrm{for}~\theta>0~\mathrm{sufficiently~small}.
	\]
	On the other hand, for all $t\in [\lambda T,T]$, it holds
	\[
	\begin{split}
		p^{\theta, \lambda, {\bf h}}(t)-\theta h
		&~\preceq~\theta\cdot\bigg[-\left({2\lambda \over (1-\lambda)}-\sqrt{\theta}\right) {\bf e}+{\lambda (1-\lambda)T^2\over 8}D^2V(0)({\bf e})\\
		&\qquad\qquad\qquad\qquad\qquad\qquad\qquad+{2(1+\sqrt{\theta})\over (1-\lambda)^2 T^2}\cdot (1,\cdots, 1)^{\dagger}+{E(t,\theta,{\bf h})\over \theta}\bigg].
	\end{split}
	\]
	By choosing $\lambda=1-T^{-4/3}$ and using (\ref{conV2}), (\ref{conT}) we derive 
	\[
	\begin{split}
		\limsup_{\theta\to 0+}\sup_{t\in [\lambda T,T]}\left({p^{\theta,\lambda,{\bf h}}(t)-\theta h\over \theta}\right)&\preceq~-{2\lambda \over 1-\lambda}{\bf e}+{\lambda (1-\lambda) T^2\over 8}D^2V(0){\bf e}+{2\over (1-\lambda)^2T^2 }\cdot (1,\cdots, 1)^{\dagger}\\
		&\preceq~\left(-{2\lambda \gamma_{\min}\over 1-\lambda}+{\lambda (1-\lambda) T^2\gamma^2_{\max}\over 8}+{2\over (1-\lambda)^2T^2 }\right)\cdot (1,\cdots, 1)^{\dagger}\\
		&\preceq~{T^{2/3}\over 8}\cdot \left(-16\lambda T^{2/3}\gamma_{\min}+\lambda \gamma^2_{\max}+16\right)\cdot (1,\cdots, 1)^{\dagger}\\
		&\preceq~{T^{2/3}\over 8}\cdot\left(-8\lambda T^{2/3}\gamma_{\min}+16\right)\cdot (1,\cdots, 1)^{\dagger}~\prec~0.
	\end{split}
	\]
	In this case, for $\theta>0$ sufficiently small, one has 
	$$
	p^{\theta,\lambda,{\bf h}}(t)~\preceq~{\bf h}\qquad\forall t\in [\lambda T, T],
	$$
	and  the map $t \mapsto \big(x^{ \theta, \lambda, {\bf h} }(t),p^{\theta, \lambda, {\bf h} }(t)\big)$ is a negative super-solution of (\ref{ex3}).
	\v

	\n {\bf 4.} To study the stability of solution $x^j(t,\xi)\equiv y_j(t)$ for $j=1,2$,  we first compute $D\Phi(y_j)({\bf b})$ for a given ${\bf b}\in \mathcal{C}([0,T],\R^d)$.  From  Step 1, the map $b\mapsto \Phi(b)$ takes the form
	\[
	\Phi(b)~=~x^b(t)\qquad\forall t\in [0,T],
	\]
	with $(x^b,p^b)$ being the unique solution of (\ref{ex2})-(\ref{bd2}). For every $\ve>0$, let $(x_{\ve},p_{\ve})$ be  $(x^b,p^b)$ with $b=y_j+\ve {\bf b}$. By the linearization, it holds $(x_\ve,p_{\ve})=(y_j,p_j)+\ve\cdot ({\bf x}_{\bf b},{\bf p}_{\bf b})+o(\ve)$ where $({\bf x}_{\bf b},{\bf p}_{\bf b})$ is the solution of the linear system 
	\[
	\begin{bmatrix}
		\dot{{\bf x}}(t)\\[2mm]
		\dot{{\bf p}}(t)
	\end{bmatrix}
	~=~
	\begin{bmatrix}
		0&-1\\[2mm]
		-2\kappa \mathbb{I}_d-D^2V(y_j(t))&0
	\end{bmatrix}
	\begin{bmatrix}
		{\bf x}(t)\\[2mm]
		{\bf p}(t)
	\end{bmatrix}
	+
	2\kappa {\bf b}(t)\begin{bmatrix}
		0\\[2mm]
		1
	\end{bmatrix},\qquad~~
	\begin{bmatrix}
		{\bf x}(0)\\[2mm]
		{\bf p}(T)
	\end{bmatrix}~=~\begin{bmatrix}
		0\\[2mm]
		0
	\end{bmatrix}.
	\]
	Hence,  we obtain an expression for the differential $D\Phi(y_j)({\bf b})={\bf x}_{\bf b}$ with ${\bf x}_b$ solving  the second order ODE
	\bel{se-ODE}
	\ddot{\bf x}(t)-\big(2\kappa\mathbb{I}_d+D^2V(y_j)\big) {\bf x}(t)+2\kappa {\bf b}(t)~=~0,\qquad  \dot{\bf x}(T)~=~ {\bf x}(0)~=~0.
	\eeq
	
	Assume that   $(\lambda, {\bf b}_{\lambda})$ is a pair of eigenvalue and eigenfunction of $D\Phi(y_j)$. Then ${\bf x}_{{\bf b}_{\lambda}}=D\Phi(y_j)({\bf x}_{{\bf b}_{\lambda}})=\lambda {\bf x}_{{\bf b}_{\lambda}}$
	solves the equation
	\bel{ei-va}
	\ddot{\bf x}(t)-\left({2\kappa(\lambda-1)\over \lambda}\cdot\mathbb{I}_d+D^2V(y_j(t))\right) {\bf x}(t)~=~0,\qquad \dot{\bf x}(T)~=~ {\bf x}(0)~=~0.
	\eeq
	For every $i\in \{1,\cdots, d\}$, let $(\gamma_i^n,z_i^n)$ be all pairs of eigenvalues and eigenfunctions for the two-point boundary problem
	\bel{eii}
	\ds \ddot{z}(t)-\left(\gamma+V_{x_ix_i}(y_j(t))\right) z(t)~=~0,\qquad\dot{z}(T)~=~ z(0)~=~0.
	\eeq
	Calling $\{{\bf e}_1,\cdot, {\bf e}_1,\cdot\}$  the standard basis of $\R^d$,  the  function  ${\bf z}^n_i\doteq z^n_i\cdot {\bf e}_i$ solves  (\ref{ei-va}) for $\lambda=\lambda_n^i$ such that $\ds {2\kappa(\lambda_n^i-1)\over \lambda_n^i}~=~\gamma_n^i$. Moreover,  by Sturm–Liouville theory, we have that the sequence $\{{\bf z}_i^n\}_{n\geq 1,i\in \{1,\cdots, d\}}$ forms a complete orthogonal basis of subspace of ${\bf L}^2\big([0,T],\R^d\big)$ functions that satisfy the boundary condition in (\ref{se-ODE}). Hence, to show that both solutions $x^1(t,\xi)$ and $x^2(t,\xi)$ are stable, we shall verify 
	\bel{ke1}
	0~<~\lambda^n_i~\leq~1\quad\forall i\in \{1,\cdots, d\}, n\geq 1.
	\eeq
	Multiplying (\ref{eii}) by $z(t)$, we get 
	\[
	{d\over dt} \big(z_i^n(t)\dot{z}_i^n(t)\big)~=~\left(\gamma+V_{x_ix_i}(x^j(t))\right) z^n_i(t)+ \dot{z}^n_i(t),
	\]
	and the boundary condition of (\ref{eii}) and (\ref{con-V}) yield
	\bel{11}
	\begin{split}
		0&~=~z_{i}^n(T) \dot{z}_i^n(T)-{z}_i^n(0)\cdot \dot{z}_i^n(0)~=~\int_{0}^{T}\left(\gamma+V_{x_ix_i}(x^j(t))\right)\big[z_i^n\big]^2(t)+ \big[\dot{z}_i^n\big]^2(t)dt\\
		&~\geq~ \int_{0}^{T}\left(\gamma+V_{x_ix_i}(x^j(t))\right)\cdot \big[z_i^n\big]^2(t)dt~\geq~\gamma\cdot \int_{0}^{T}\big[z_i^n\big]^2(t)dt.
	\end{split}
	\eeq
	Since $z_i^n$ is a nonzero solutution, we obtain 
	\[
	{2\kappa(\lambda_n^i-1)\over \lambda_n^i}~=~\gamma_n^i~\geq~0\quad \forall i\in \{1,\cdots, d\}, n\geq 1,
	\]
	and this  yields (\ref{ke1}) and the negative  solution $x^2(t,\xi)=y_2(t)$ is stable.
	\medskip

	In the case $j=1$, using  (\ref{conV2}) in (\ref{11}), we derive 
	\[
	0~\geq~ \int_{0}^{T}\left(\gamma+V_{x_ix_i}(0)\right)\cdot \big[z_i^n\big]^2(t)dt~\geq~(\gamma+\gamma_{\min})\cdot \int_{0}^{T}\big[z_i^n\big]^2(t)dt,
	\]
	which implies that 
	\[
	{2\kappa(\lambda_n^i-1)\over \lambda_n^i}~=~\gamma_n^i~<~-\gamma_{\min}.
	\]
	Solving the above inequality, we get 
	\[
	0~<~\lambda_n^i~<~{2\kappa\over 2\kappa+\gamma_{\min}}~<~1 \quad \forall i\in \{1,\cdots, d\}, n\geq 1,
	\]
	and the zero  solution $x^1(t,\xi)=y_1(t)=0$ is asymptotically stable.  
	\endproof

	\v

	{\bf Acknowledgements.} This research was partially supported by NSF-DMS 2510261 (L. Bociu) and NSF-DMS 2154201 (K.T. Nguyen). K.T. Nguyen would also like to warmly thank Professor A. Bressan for suggesting the problem.
	

\end{document}